\documentclass{ws-procs10x7square}

\def\BbC{{\rm C\!\!\!I\,}}
\newcommand{\BbR}{{\rm I\!R}}

\def\un { \rm 1\kern-5.8ptI}

\newcommand{\beq}{\begin{eqnarray}}
\newcommand{\eeq}{\end{eqnarray}}
\newcommand{\bq}{\begin{equation}}
\newcommand{\eq}{\end{equation}}
\newcommand{\Sum}{\displaystyle \sum}
\newcommand{\Int}{\displaystyle \int}
\newcommand{\Frac}{\displaystyle \frac}
\def \1{{\rm 1\kern-4pt 1}}

\newcommand{\Min}{\displaystyle \min}

\begin{document}


\title{Dirac-Fock models for atoms and molecules and related topics}

\author{Maria J. Esteban and Eric S\'er\'e}
\address{Ceremade (UMR CNRS 7534), Universit\'e Paris IX-Dauphine, 75775 Paris C\'edex~16, France. esteban, sere@ceremade.dauphine.fr}


\maketitle

\abstracts{
}


\section{Introduction}

{}\ 

 Relativistic effects are important in the electrons' dynamics and bound state energies in heavy atoms and molecules. When the nucleii involved are heavily charged, the velocities of the electrons of the inner layers are quite large, and so nonrelativistic modelling will lead to important errors. The usual strategies  to address this issue are : either use nonrelativistic models together with relativistic corrections, or use  relativistic models based on the Dirac operator.

\smallskip
The free Dirac operator is the simplest first order (in time and space) constant coefficient operator which is invariant under the action of the Lorentz group, {\sl i.e.} its is compatible with Relativity Theory. The free Dirac operator can be written as: 
$$i\partial_t+H_c\,,$$
where
\[{{ H}_c = -i\,c\,\hbar \,\alpha \cdot \nabla +
mc^2\beta}\;,\quad\mbox{with}\; {\alpha_1,\; \alpha_2,\;
\alpha_3,\; \beta} \in {\mathcal M}_{4 \times 4} (\BbC)\]
\[{
\beta =
\left ( \begin{array}{cc} \un & 0 \\
0 & - \un \\ \end{array} \right ) \; , \ \alpha_i= \left (
\begin{array}{cc} 0 &
\sigma_i \\ \sigma_i & 0 \\ \end{array} \right )}\]
and the matrices $\,\sigma_k\,$ are the Pauli matrices:
\[{\sigma _1=\left( \begin{matrix}0 & 1
\\ 1 & 0 \end{matrix}\right) \; ,\quad \sigma_2=\left(\begin{matrix} 0 & -i \\
i & 0 \end{matrix}\right) \; ,\quad \sigma_3=\left( 
\begin{matrix} 1 & 0\\ 0 &-1\end{matrix}\right)}\] 

Two of the main properties of $H_c$ are that its square is a shifted Schr\"odinger operator:
\bq\label{carre}H_c^{^2}=-c^2\hbar^2\Delta +m^2c^4\,,\eq
and that its spectrum is unbounded below and above:
\bq\label{unbdd} \sigma(H_c)=(-\infty, -mc^2] \cup [mc^2, +\infty)\,,\eq
and therefore has a gap, $(-mc^2, mc^2)$.

\smallskip
Looking for bound states for an electron in an external potential $V$ of the form
$$ \tilde\psi (t,x)\; =\; e^{-i\lambda t}\, \psi (x)
 $$
is equivalent to looking for eigenvalues  $\,\lambda\,$ of the operator $\,H_c+V$ associated with the eigenfunctions $\,\psi$. In the case of $N$ electrons, by analogy with the nonrelativistic case, the Hamiltonian should be written as : $\,H_{c,x_1}+\cdots+,H_{c,x_N}+V\,$. But as Brown and Ravenhall remarked in \cite{Brown-Ravenhall-51}, for potentials $V$ not too singular,  the spectrum of this Hamiltonian is the whole real line. So, no isolated eigenvalues can correspond to bound state energies. But by analogy to the nonrelativistic case, one could still try to find eigenfunctions of $\,H_{c,x_1}+\cdots+,H_{c,x_N}+V\,$ as Slater determinants, that is, in the form:
\bq\label{slater}  {\psi:=\frac{1}{\sqrt{N!}}\det(\psi_i(x_j))\,,\quad i,j=1,\dots, N\,,}\eq
the functions $\,\psi_j$ being mutually orthogonal and of norm $1$ in the space $\,L^2(\BbR^3,\BbC)$ (the constant $\,\frac1{\sqrt(N!)}\,$ appearing in the above expression  just for the sake of normalization).
The reason to consider Slater determinants is the Pauli principle which states that one cannot find two electrons  in the same state. Hence, the eigenfunctions should be antisymmetric and the simplest antisymmetric function is the determinant.

  \smallskip
The Dirac-Fock equations can be found in the above manner. In that sense, they are the relativistic counterpart of the nonrelativistic   Hartree-Fock equations which are derived in the same way from the Hamiltonian
$\,-\Delta_{x_1}-\cdots -\Delta_{x_N}+V\,$.  Despite their non very physical derivation, the Dirac-Fock model is widely used in computational atomic physics and quantum chemistry to study atoms and molecules involving heavy nucleii. These equations were first introduced by Swirles in \cite{Swirles-35}. One can find many articles in the Physics litterature about the Dirac-Fock equations (or more complicated models, like the multi-configuration Dirac-Fock model  in the case of atoms and molecules with open shells) : see for instance  
\cite{Desclaux-93, Desclaux-Mayers-O'Brien-71, Dyall-Grant-Wilson-84, Ellis-Painter-70, Froese Fischer-77, Grant-82, Grant-89, Grant-87, Grant-Quiney-88, Heully-Lindgren-Lindroth-Lundqvist-
Maartensson-Pendrill-86A, Indelicato-Desclaux-93, Indelicato-95, Johnson-Blundell-Sapirstein-88, Johnson-Lin-76, Kullie-Dusterhoft-Kolb-99, Kullie-Kolb-01, Paturel-00, Quiney-Grant-Wilson-87}.

  \smallskip
An interesting question is why these equations have bound state solutions (see below)  and give very satisfactory numerical results.  Why are these equations so well fitted to describe the stable electronic configuration in atoms and molecules if the model is not physically ``well posed"? In Section 3 we will discuss a possible link between  the Dirac-Fock equations and models of Quantum field theory which could maybe explain these phenomena.

  \smallskip

 A function {$\Psi\in \Sigma:=\{\Psi=(\psi_1, \dots, \psi_N)\,;\;( \psi_\ell, \psi_k)_{_{L^2}} \ = \ \delta_{_{k\ell}}\} $} is said to satisfy the Dirac-Fock equations if
   for $\,k=1,\dots N$, 
\bq\label{DF}  \qquad\bar H_{_{c,\Psi}} \, \psi_{_k} \ = \  \lambda_{_{k}} \
\psi_{k}  \ , \quad \lambda_k\in (-mc^2, mc^2)\,,\eq
where  the mean-field operator  $\,\bar H_{_{c,\Psi}}\,$ is defined by
\bq\label{meanfield}  \bar H_{_{c,\Psi}}\, \psi_k = \left(H_c+V\right) \psi_k
+  (\rho \ast \frac{1}{|x|}) \psi_k -
 \int_{\BbR^3}\frac{ R(x,y) \psi_k(y)} {|x-y|}\,dx\,dy \,,\eq
  {$\rho$} being  the (scalar)   {electronic density}
and the  $4 \times 4$  complex matrix   {$R$},    {the exchange matrix} which  comes from the antisymmetry
of the Slater determinant :
$$
  {\rho (x) \ = \ \displaystyle{\sum^N_{\ell = 1}} \ \Bigl( \psi_{_\ell} (x),
 \psi_{_\ell} (x) \Bigr)_{\BbC^4}}\;,\quad
  {R(x,y) \ = \ \displaystyle{\sum^N_{\ell = 1 }} \ \psi_{_\ell} (x) \otimes 
\psi^\ast_{_\ell} (y) \,.}$$

  \smallskip
Formally, the solutions of DF equations are the stationary points of the ``energy functional"
$$  {{\mathcal E_{c}} (\Psi) = \displaystyle{\sum^N_{k = 1}} 
\Bigl((H_c+V)\psi_{_k}, 
\psi_{_k} \Bigr)_{_{L^2}}+  \frac{1}{2} \ \displaystyle  {\iint_{_{\BbR^3 
\times \BbR^3}}}\frac{
 \rho(x) \rho(y) -  |R(x,y)|^2 }{|x-y]}\, dx \,dy \,,
}$$

\noindent in the set $$\,\Sigma:=\{\Psi=(\psi_1, \dots, \psi_N)\,;\;( \psi_\ell, \psi_k
)_{_{L^2}} \ = \ \delta_{_{k\ell}}\} \,.$$

\noindent{\bf Remark.} The infimum of   $\,{\mathcal E}_c \,$   over the set $S$ is equal to $\,-\infty\,$
So, the solutions of DF will have to be some kind of saddle points for $\,{\mathcal E}_c \,$.

\section{The Dirac-Fock equations and their nonrelativistic limit}

In an atomic or molecular model, the external potential can be aproximately described by 
$$V=-Z\mu*\frac 1{|x|}\,,$$
where the  measure $\,Z\mu\,$ represents the nuclear charge distribution, and in the particular case of a molecule with $K$ point nucleii  with charge $\,Z_i$ and situated at the points 
$\,{\bar x}_i, i=1,\dots K\,$,  the measure $Z\mu$ takes the form $\,Z\mu=\Sum_i Z_i\delta_{{\bar x}_i}$. So, $Z$ is the total nuclear charge.

  \smallskip
As already said in the Introduction, finding solutions of (\ref{DF}) is then reduced to finding critical points of the functional $\,{\mathcal E_{c}}$ on the set $\Sigma$. Once again, the unboundedness (from above and below) of the spectrum of the free Dirac operator makes the functional $\,{\mathcal E_{c}}\,$ totally indefinite. This together with the {\sl a priori} lack of compactness of the problem ($\,\BbR^3\,$ is unbounded) and the fact that we have to work on a manifold $\Sigma$ and not in the whole functional space, makes the variational problem difficult. A minimization procedure is once again out of the question and another method has to be found. In \cite{Esteban-Sere-99} the authors defined a penalized variational problem which could be solved by first maximizing on some part of the spinor functions $\,\psi_i$ and then defining a more standard min-max argument for a reduced functional. 

  \smallskip
The theorem proved in \cite{Esteban-Sere-99} states the following :
\begin{theorem}\label{thmDF} {\rm(\cite{Esteban-Sere-99}) } With the above notations, assume that $N$ and $Z$ are two positive integers satisfying $\, 
 \max (Z, 3N-1) < \frac{2c^2}{\pi/2+2/\pi}$ and $\,N<Z+1$. Then, there exists an infinite  sequence $(\Psi^{c,j})_{_{j \geq 1}}$ of critical points of the
DF functional ${\mathcal E_{c}}$ on $\Sigma$.
 The function vectors  $\,(\psi^{c,j}_1, \dots \psi^{c,j}_N)\,$ belong to $\, \Sigma\, $  
and they are strong solutions, in
$H^{1/2}(\BbR^3, \BbC^{\,4}) \ \cap \ \bigcap_{1\leq q < 3/2} W^{1,q}(\BbR^3,
\BbC^{\,4})$, of the Dirac-Fock equations, that is, for all $\,j\geq 1$, 
\begin{equation} \bar{H}_{_{\Phi^{c,j}}} \ \psi^{c,j}_{_k} \ = \ \epsilon^{c,j}_{_k} \
\psi^{c,j}_{_k} \ , \quad 1 \leq k \leq N \,,\end{equation} \begin{equation} 0 \ < \
\epsilon^{c,j}_{_1} \ \leq ... \leq \ \epsilon^{c,j}_{_N} \ < \ mc^2 \ . \end{equation} Moreover,
\begin{equation} 0 \ < \ {\mathcal E_{c}} (\Phi^{c,j}) \ < \ Nmc^2 \ , \end{equation}
\begin{equation} \lim_{j \rightarrow \infty}  {\mathcal E_{c}} (\Phi^{c,j}) \ = \ Nmc^2 \ . \end{equation} 
\end{theorem}

  \smallskip

\begin{remark} Since $\mu$ is arbitrary, the assumptions of the above theorem contain the case of
point-like nuclei as well as more realistic nuclear potentials  of the
form
 $- \displaystyle{ \sum_i\rho_i(x) *\frac{1 }{ \vert x \vert}}\,$, 
where
$\,\rho_i\in L^\infty\cap L^1\,,\; \rho_i\geq 0\,, \, \displaystyle{\sum_i
\int_{\BbR^3}\rho_i=Z}\,.$
\end{remark}

\begin{remark}  In our units, the above conditions become $$Z \leq 124 \ , \ N \leq 41 \ , \quad N \leq Z \ .$$ 
\end{remark}

  \smallskip
The proof of the above theorem is done by defining a sequence of variational problems corresponding to increasing 'topological complexity' or Morse index.
By the application of our 'first variational argument' we find solutions $\,\Psi^{c,1}$,  that will play an important role (see below), since for $c$ large they will be actually 'electronic ground states'. 

  \smallskip
Soon after this result was published, E. Paturel used an alternative method  to prove the same result  but without the unnatural condition 
\bq\label{unnatural}\,  \max (Z, 3N-1) < \frac{2c^2}{\pi/2+2/\pi}\,,\eq (see \cite{Paturel-00}). Note that the assumption $N<Z+1$ is already present in all the existing results for the nonrelativistic Hartree-Fock equations. Also, the conditions
$\,Z$ and $\,N$ less than $\,2c^2/(2\pi+\pi/2)\,$ are not that unnatural, since already in the linear case such a condition is necessary to use Hardy-like inequalities ensuring the existence of a gap of the spectrum of $H_c+V$ around $0$ (see for instance \cite{Tix-97, Tix-98} and \cite{Dolbeault-Esteban-Sere-00A}). But condition (\ref{unnatural}) was clearly related to the particular method of proof in \cite{Esteban-Sere-99}. Paturel's theorem gets rid of (\ref{unnatural}) and proves the same as Theorem \ref{thmDF},  but under the sole assumptions
$$\,Z< \frac{2c^2}{\pi/2+2/\pi},\quad \,N< \frac{2c^2}{\pi/2+2/\pi},\quad  N<Z+1\,.$$

  \smallskip
Natural questions which arise from the above results are the following:

\smallskip
- Is there a notion of ground-state in this model?

\smallskip
- Is there a link between these solutions and the solutions of the usual Hartree-Fock solutions?

\smallskip
- Can we find the solutions in a simple way? By a ``simple" variational argument?

\medskip
To answer the above questions, it is useful to investigate
 the nonrelativistic limit of the Dirac-Fock equations.
Of course, their formal limite are the well-known Hartree-Fock equations, which are identical to the Dirac-Fock ones, but replacing the Dirac operator $\,H_c\,$ with the Schr\"odinger operator $\,-\frac{\Delta}2$. In \cite{Esteban-Sere-01}, this  was proved rigorously, and as we see below, this result has been of importance to better understand the variational structure of the Dirac-Fock problem and in particular to obtain a good definition of an electronic ground-state energy, which is {\sl a priori} not clear because of the unboundedness of the Dirac-Fock energy. 

  \smallskip
The following is proved in  \cite{Esteban-Sere-01}.

\begin{theorem}\label{A}  Let $N < Z+1$.
Consider a sequence of numbers $c_n\to +\infty$ and a sequence
$\{\Psi^n\}_n$ of solutions of  (\ref{DF}), i.e.
$\Psi^n=(\psi_1^n,\cdots,\psi^n_N)$, each $\psi_k^n$ being
in $H^{1/2}(\BbR^3,\BbC^4)$, with $\Int_{\BbR^3}
{\psi_k^n}^\ast
\psi_l^n\,dx = \delta_{kl}$ and $\bar H_{_{c_n,\Psi^n}}\psi^n_k
=\varepsilon^n_k\psi^n_k$ . Assume also  that the
multipliers $\varepsilon^n_k$  \break ($\,k=1,\dots, N \,$) satisfy
$$0 < c^2-\mu_j < \varepsilon^{^{c,j}}_1 \leq ... \leq
\varepsilon^{^{c,j}}_N <
c^2-m_j , \;\hbox{ with $\;\mu_j > m_j > 0$ independent of $\,n$}\,. $$
 Then for $n$ large enough, $\psi_k^n$ is
in $H^{1}(\BbR^3,\BbC^4)$  and there exists a solution of the Hartree-Fock equations,
$\bar \Phi = \left( \bar
\varphi_1, \cdots , \bar \varphi_N \right) $, with negative
multipliers,
$\bar{\lambda}_1, ...,
\bar{\lambda}_N\,$, such that, after extraction of a subsequence,

\bq \lambda^{n}_k:=\,\varepsilon^n_k - (c_n)^2\quad
\longrightarrow_{_{\hskip-7mm n\to +\infty}}\;
\bar{\lambda}_k
\ ,
\quad k = 1, ...,N \ ,
     \label{(4)} \eq

\beq \psi_k^{n} = \left({\begin{array}{c}
\varphi_k^{n}
\\ \chi_k^{n}\end{array}}\right) \longrightarrow_{_{_{\hspace{-6mm} n \rightarrow + \infty
}}} \left( {\begin{array}{c}\bar
\varphi_k
\\ 0\end{array}}\right)
   \ \; {\rm in \ } \ H^1 (\BbR^3,\BbC^{\,2})\times H^1 (\BbR^3,\BbC^{\,2}),
\label{(5)}
\eeq

\beq \left\Vert\chi_k^{n}  +\frac {i}{2c_n}(\sigma\cdot
\nabla)\varphi_k^{n}\right\Vert_{_{L^2(\BbR^3,\BbC^{\,2})}}=O(1/(c_n)^3 ),
\label{(5bis)}
\eeq
and
\bq {\mathcal E}_{c_n}(\Psi^n)-Nc_n^2\quad
\longrightarrow_{_{_{\hspace{-7mm} n
\rightarrow + \infty }}} \quad {\mathcal E}_{HF}(\bar
\Phi). \label{(5ter)}
\eq
\end{theorem}

  \smallskip
Actually the above theorem can be made more precise :  the sequence of the solutions $\,\Psi^{c,1}\,$ given by Theorem \ref{thmDF} are shown to converge towards ground state solutions of the Hartree-Fock problem when $c$ goes to $+\infty$. This gives a ``particular" status to that solution of the Dirac-Fock equations.  It does not minimize the Dirac-Fock energy among all functions in the set $\,\Sigma$, since $\,\inf_\Sigma {\mathcal{E}_{c}}=-\infty$, but when taking the nonrelativistic limit, they approach those solutions of the Hartree-Fock equations which minimize the corresponding Hartree-Fock energy. Actually, there is more to it. In \cite{Esteban-Sere-01} the following theorem was proved :

\begin{theorem}\label{M} {\rm(\cite{Esteban-Sere-01})}  Fix $N,Z$ with $N<Z+1$ and take $c$ sufficiently
large.
Then $\Psi^{c, 1}$ is a solution of the following minimization problem:
\bq\label{minpro1}\inf\{{\mathcal E}_c(\Psi)\,;\; {\rm Gram_{_{L^2}}}\Psi=
\un_{_N} \,,\; \Lambda^{^-}_{\Psi}\,\Psi=0 \ \}\eq where
$\Lambda^{^-}_{\Psi}=\chi_{(-\infty,
0)}(\bar H_{c,\Psi})$ is the negative spectral projector of the
operator
$\bar
H_{c,\Psi}$, and $\Lambda^{^-}_{\Psi}\,\Psi:=(\Lambda^{^-}_{\Psi}\,
\psi_1,\cdots ,\Lambda^{^-}_{\Psi}\,\psi_N)\;.$
\end{theorem}
Actually,  in \cite{Esteban-Sere-01} we call {\sl electronic configurations} those functions $\,\Psi\,$ in $\,\Sigma\,$ which satisfy $\,  \Lambda^{^-}_{\Psi}\,\Psi=0\,$ (note that all solutions of the Dirac-Fock equations with positive eigenvalues $\ \varepsilon_k\,$ are electronic configurations by definitio). Thus, the above result shows that the Dirac-Fock energy is bounded from below in the set of all electronic configurations, and $\,\Psi^{c, 1}$ is a minimizer for it in that set. 

  \smallskip
The above theorems answer some of the questions made above. Another question that we asked was that of finding a non very complicated variational argument to find the solutions exhibited in Theorem \ref{thmDF}, since the proofs in that theorem are very
involved technically and very difficult to implement in actual computations.
In order to do this, in \cite{Esteban-Sere-01} we introduced the notion of projector ``$\varepsilon$-close
to $\Lambda^+_c$'', where $\Lambda^+_c = \Frac{1}{2} \big | H_c \big
|^{-1} \bigl( H_c +
    \big | H_c \big | \Bigr)$ is the positive free-energy spectral projector.

\begin{definition} \label{P} : Let $P^+$ be an orthogonal projector
in $L^2 ( \BbR^3, \BbC^4 )$, whose restriction to
$H^{\frac{1}{2}} ( \BbR^3,
\BbC^4 )$ is a bounded operator on $H^{\frac{1}{2}} (
\BbR^3,\BbC^4 )$ .

\noindent Given $\varepsilon > 0$, $P^+$ is said to be $\varepsilon$-close  to
$\Lambda^+_c$ if and only if, for all
$\psi \in
H^{\frac{1}{2}} (\BbR^3, \BbC^4)$,
$$
\Big \Vert \Bigl(-c^2 \Delta + c^4 \Bigr)^{\frac{1}{4}} \Bigl(P^+ -
\Lambda^+_c \Bigr) \psi \Big \Vert _{L^2 (\BbR^3, \BbC^4)}
    \leq  \varepsilon \  \Big \Vert \Bigl(-c^2 \Delta + c^4
\Bigr)^{\frac{1}{4}}  \psi \Big \Vert _{L^2 (\BbR^3, \BbC^4)} \ .$$
\end{definition}

  \smallskip
An obvious example of projector $\varepsilon$-close to $\Lambda^+_c$ is
$\Lambda^+_c$ itself. Other interesting examples are the mean-field  operators $\,\bar H_{c,\Psi}$ for $\,c\,$ large enough. Let us now
give a min-max principle
associated to $P^+$ :
\begin{theorem} \label{minma} Fix $N, Z$ with $N < Z+1$. Take $c >
0$ large enough, and $P^+$ a projector $\varepsilon$-close  to
$\Lambda^+_c$, for $\varepsilon > 0$
small enough. Let $P^- = \un_{_{L^2}} - P^+$, and define
$$ E(P^+) : = \inf_{{\scriptstyle \Phi^+ \in (P^+ H^{\frac{1}{2}} ) ^{^N}}
\atop \scriptstyle {\rm Gram_{_{L^2}}} \Phi^+ = \un_N} \quad
\sup_{{\scriptstyle\Psi \in (P^- H^{\frac{1}{2}} \oplus \,{\rm Span}
(\Phi^+))^N \atop \scriptstyle {\rm Gram_{_{L^2}}} \Psi = \un_{_N} }}
\quad {\mathcal E}_c (\Psi)\ . $$
Then $E(P^+)$ does not depend on $P^+$ and $\,{\mathcal E}_c (\Psi^{c,1} )
= E(P^+)$.

\end{theorem}

This result shows that when we want to find the electronic
ground-state of the Dirac-Fock problem by a ``simple min-max", we have a large choice of projectors   in the limit $\,c\to +\infty$. Other characterizations
of the ground-state energy for $c$ large can be found in \cite{Esteban-Sere-01-2}.

 \section{The Dirac-Fock equations and a Hartree-Fock model in Fock's space}

\smallskip
Atomic and molecular models based on Q.E.D. seem to be better motivated than the Dirac-Fock equations.
In a Q.E.D. formulation one should be able to write an energy functional which is bounded below in a well chosen subset of the Fock space. In such a model one cannot fix the number of particles (electrons) but only the charge, since one has to admit the possible presence of electron-positron pairs, which do not change the charge, but change the number of particles. In all models based on Q.E.D tehe are some choices to be made, for instance,  a notion of electron, that is, a positive energy space (and its corresponding projector). Also, the kind of configurations that are allowed in the model has to be made precise, in order to see if multiconfiguration is allowed or not.

  \smallskip
In \cite{Bach-Barbaroux-Helffer-Siedentop-98, Bach-Barbaroux-Helffer-Siedentop-99} such a model has been
proposed and the case of no electrons completely treated (the empty set). In \cite{Barbaroux-Farkas-Helffer-Siedentop-xx} the same kind of model has been used to treat the case of positive charges ($N\geq 0$). More concretely, in these models
one chooses   {$\,P^+ \in {\mathcal P}$,  ``a positive projector" (choice of electronic space)}, and
$$
  {S^{ N, P^+}  =  \{\gamma\in {S}_1 (L^2)\, ,\
\gamma=\gamma^*,\, -(1-P^+) \leq \gamma \leq P^+, \;P^+\,\gamma\,(1-P^+) =0,
\; \mbox{ tr}\,\gamma\leq N \}\,,}$$
where $S_1(L^2)$ is the trace class operators Hilbert space.

\smallskip
Let us now define
$$ { {\mathcal F}^c_{\alpha}(\gamma) = {\rm tr}\left((H_c +V- mc^2)\gamma\right) +
\frac{{1}}{2} \int \frac{\rho_\gamma({\bf x})\rho_\gamma({\bf
y})-|\gamma(x,y)|^2}{| {\bf x} - {\bf y}|}\,{\rm d}{\bf x}\, {\rm d}{\bf y}}\,,
$$ 
\bq\label{min} I(c,N, P^+)= \Min_{_{ \gamma\in S^{N, P^+}}}  \; {\mathcal F}_c   (\gamma)\,.\eq

\smallskip
 In a recent result, Barbaroux, Farkas, Helffer and Siedentop have proved that 
 when one chooses the class of projectors
 $$ {{\mathcal P} := \{ P^+_{\Psi} ,\ P^+_{\Psi} =
 \chi_{[0, \infty)} (\bar H_{c, \Psi})\;, \quad \Psi\in \Sigma\},}$$
the minimization problem (\ref{min}) has a solution. Moreover,  under some (reasonable) conditions,  {the minimizer of $\,{\mathcal F}_c\,$ in $\,S^{N, P^+}\,$  is { no-pair (purely electronic)}, that is, there exists functions $\,\bar\psi_1, \dots, \bar\psi_N\,$ such that : $\displaystyle{\,\gamma_{_{min}}
= \sum_{k=1}^N \langle\bar\psi_k\, , .\rangle \bar\psi_k
\quad \mbox { and}} \qquad$} 
\bq\label{projDF} {P^+}\,(\bar H_{c,\bar\Psi})\, {P^+}\,\bar\psi_k=\lambda_k\,\bar\psi_k\,,\qquad \lambda_k\in (0, mc^2)\,,\eq
 that is, the $\,\bar \psi_k$'s are solutions to the projected Dirac-Fock equations (projected in the image of the chosen projector $P^+$, of course).
 
   \smallskip
 Now, in  \cite{Mittleman-81}, Mittleman proved (somehow formally, but the main ingredients of a proof  are present in the paper) that any stationary solution (stationary in $P^+$ and $\gamma$) of a related QED functional would be a solution of the self-consistent Dirac-Fock equations. Hence, when we take the energy $\,{\mathcal F}_c\,$ to depend both on $\,P^+\,$ and  $\,\gamma\,$, any stationary point of this new functional, $\,(\bar\Psi, \bar P^+)\,$,  should be a solution of (\ref{projDF}) with
 $\bar P^+=\chi_{_{(0, +\infty)}}(\bar H_{c,\bar\Psi})$.
 
   \smallskip
 A natural idea would be to find stationary points in  $\,P^+\,$ and  $\,\gamma\,$ by considering the max-min problem:
\bq\label{maxmin} {\sup_{P^+\in {\mathcal P}}\;\; \inf_{\gamma\in S^{N, P^+}} {\mathcal F}_c(\gamma)}\,.\eq
 
   \smallskip
This problem was shown to have a unique solution  in the  case $\,N=0\,$   in \cite{Bach-Barbaroux-Helffer-Siedentop-98, Bach-Barbaroux-Helffer-Siedentop-99}. In the case $\,N>0\,$, Barbaroux, Farkas, Helffer and Siedentop prove that indeed Mittleman's result holds in their context. It is then natural to ask whether the max-min problem defined above has a solution or not, and in the  affirmative case, whether the solution of that variational problem is a stationary point of the energy functional or not.
 
 \medskip
 In a recent  joint work with Barbaroux (\cite{Barbaroux-Esteban-Sere-xx}), we prove that the answer to the above question is {\sl yes} in some particular cases
(for instance when the ``neighboring"  linear case corresponds to a closed shell atom), but it is {\sl no} in the general case; indeed, we exhibit a case  where even if the max-min defined above were attained by some  $\,(\bar\Psi, \;\bar P^+)\;$, $\,(\bar\Psi, \;\bar P^+)\;$
 would not be a solution of the self-consistent Dirac-Fock equations, that is $\;\bar P^+\ne \chi_{(0, +\infty)}(\bar H_{c,\bar \Psi}) \,$.
  
  \smallskip
  More precisely, it is well known that the Dirac operator is invariant under the action of the Lorentz group. This means that for every spatial rotation $\,R$, there is a unique matrix $ \,U(R)\in SU_2\,$ such that
the Dirac Hamiltonian is invariant under the transformation
  $$\psi(\cdot) \mapsto A(R)\bullet\psi:=\left(\begin{matrix} U(R)&0\\ 0&U(R) \end{matrix}\right)\psi(R\,\cdot)\,,$$
  and viceversa.
  
    \smallskip
  Then, we prove in \cite{Barbaroux-Esteban-Sere-xx} the following

\begin{proposition}\label{closedshell} Assume that $\,\mu_{c,k}\,$ (resp. $N_{c,k}$) denotes the $k$-th eigenvalue of the operator $\,H_c+V$, not counted with multiplicity (resp. its multiplicity). If  $\, N=\Sum_{k=1}^I  N_{c,k}\,$ for some positive integer $I$, then
for $c$ large enough, the max-min defined in (\ref{maxmin}) is achieved by a solution of the (self-consistent) Dirac-Fock equations, that is, if we denote by  $\,(\bar\Psi, \;\bar P^+)\;$ a solution of (\ref{maxmin}), then 
$\;\bar P^+= \chi_{(0, +\infty)}(\bar H_{c,\bar \Psi}) \,$ for $c$ large enough. \end{proposition}
 
 \noindent{\bf Remark.} The main reason why the above proposition holds true is that under the assumptions of the proposition, when $c$ is large enough, for every spatial rotation $R$, for every $\,\Psi\in \Sigma$, $\Psi$ solution of the Dirac-Fock equations,  $ {\,\bar H_{c,\Psi}=\bar H_{c, \,A(R)\bullet\Psi}}$.
 
  \smallskip
In fact we conjecture that when there is no positive integer $I$ such that  $\, N=\Sum_{k=1}^I  N_{c,k}\,$ (open shell case),
the above proposition is not true. Actually, in a special case where $V$ is a ``small" Coulomb potential, $c$ is large enough and $\, N=1+\Sum_{k=1}^I  N_{c,k}\,$, $I$ being any positive integer, we have proved that Proposition \ref{closedshell} does not hold.

  \smallskip
The above results show that there is a clear link between the Dirac-Fock equations and some variational problem issued from Q.E.D.  since in some cases both methods have the same solutions.  But this is not the case for every electronic number.


\begin{thebibliography}{0}

\bibitem{Bach-Barbaroux-Helffer-Siedentop-98} V. Bach, J.M. Barbaroux, B. Helffer, H. Siedentop. {\it  Doc. Math.} {\bf 3} 353 (1998)

\bibitem{Bach-Barbaroux-Helffer-Siedentop-99} V. Bach, J.M. 
Barbaroux, W. Farkas,  B. Helffer, H.
Siedentop.  {\it  Comm. Math. Phys.} {\bf 201}(2)
445 (1999)

\bibitem{Barbaroux-Esteban-Sere-xx} J.M. Barbaroux, M.J. Esteban, E. S\'er\'e. In preparation.

\bibitem{Barbaroux-Farkas-Helffer-Siedentop-xx}  J.M. 
Barbaroux, B. Helffer, H. Siedentop. Preprint.

\bibitem{Brown-Ravenhall-51} G. E. Brown, D. E. Ravenhall. {\it Proc. Roy. Soc.
London}  {\bf A208} 552 (1951)

\bibitem{Desclaux-93} J. P. Desclaux.
 {\it Methods and Techniques in Computational
Chemistry Clementi}. E. Ed \textsl{A: Small Systems} STEF (1993)

\bibitem{Desclaux-Mayers-O'Brien-71} J. P. Desclaux, D. F. Mayers, F. O'Brien. {\it 
J. Phys. B: At. Mol. Opt. Phys} \textbf{4} 631 (1971)

\bibitem{Dolbeault-Esteban-Sere-00A} J. Dolbeault, M.J.
Esteban, E. S\'{e}r\'{e}. {\it Calc. Var. and P.D.E.} {\bf 10} 321 (2000)

\bibitem{Dyall-Grant-Wilson-84} K.G. Dyall, I.P. Grant, S. Wilson. 
{\it J. Phys. B: At. Mol. Phys}
{\bf 17} 493 (1984)

\bibitem{Ellis-Painter-70} D. E. Ellis, G. S. Painter. {\it
Phys. Rev. B}  \textbf{2(8)} 2887  (1970)

\bibitem{Esteban-Sere-99} M.J. Esteban, E. S\'{e}r\'{e}. { 
\it Comm. Math. Phys.} {\bf 203} 499 (1999)

\bibitem{Esteban-Sere-01} M.J. Esteban, E. S\'{e}r\'{e}. { 
\it  Ann. H.Poincar\'e}  {\bf 2} 941 (2001)

\bibitem{Esteban-Sere-01-2} M.J. Esteban, E. S\'{e}r\'{e}. { \it
 Contemp. Math.} {\bf 307}135 (2002)

\bibitem{Froese Fischer-77} C. Froese Fischer. 
{\it The Hartree-Fock Method for Atoms}. Wiley, 1977.

\bibitem{Grant-82} I.P. Grant. {\it Phys. Rev. A}  {\bf 25(2)} 1230 (1982)

\bibitem{Grant-89} I.P. Grant. {\it  A.I.P. Conf. Proc.} {\bf 189} Ed. P.J. Mohr, W.R. Johnson,
J.~Sucher.  235 (1989) 

\bibitem{Grant-87} I. P. Grant. {\it
Meth. Comp. Chem.} \textbf{2} 132 (1987)

\bibitem{Grant-Quiney-88} I. P. Grant, H. M. Quiney. {\it
Adv. At. Mol. Phys.} \textbf{23} 37 (1988)

\bibitem{Heully-Lindgren-Lindroth-Lundqvist-
Maartensson-Pendrill-86A} J. L. Heully, I. Lindgren, E. Lindroth, S. Lundqvist,
A. M. M{\aa}rtensson-Pendrill {\it J. Phys. B: At. Mol. Phys.} \textbf{19} 2799 (1986)


\bibitem{Indelicato-Desclaux-93} {\it Phys. Scr.}  \textbf{T46} 110 (1993)

\bibitem{Indelicato-95} P. Indelicato. Projection Operators in
{\it Phys. Rev.~A} \textbf{51(2)} 1132 (1995)

\bibitem{Johnson-Blundell-Sapirstein-88} W. R. Johnson, S. Blundell, 
J. {\it Phys. Rev. A}  \textbf{37(2)} 307 (1988)

\bibitem{Johnson-Lin-76} W. R. Johnson, C. D. Lin. {\it Phys. Rev. A} \textbf{14(2)} 565 (1976)

\bibitem{Kullie-Dusterhoft-Kolb-99} O. Kullie, C. D\"usterh\"oft, D. Kolb. 
{\it Chem. Phys. Lett.} {\bf 307} 259 (1999)

 \bibitem{Kullie-Kolb-01}  O. Kullie{\it  Eur. Phys. J. D} {\bf 17 (2)} 167  (2001)    
 
 \bibitem{Mittleman-81} M.H. Mittleman. Theory of relativistic effects on atoms:
Configuration-space Hamiltonian.
Phys. Rev. A {\bf 24(3) } (1981), p. 1167-1175.  
                                  
\bibitem{Paturel-00} E. Paturel. {\it Ann. H.Poincar\'e} {\bf 1} 1123 (2000)

\bibitem{Quiney-Grant-Wilson-87} H.M. Quiney, I.P. Grant, S. Wilson. 
{\it  J. Phys. B: At. Mol. Phys.} {\bf 20 } 1413 (1987)

\bibitem{Swirles-35} B. Swirles.{\it
Proc. Roy. Soc. A}  {\bf 152 } 625 (1935)

\bibitem{Tix-97} C. Tix. {\it
Phys. Lett. B} {\bf 405 } 293  (1997)

\bibitem{Tix-98} C. Tix. {\it  Bull. London Math. Soc.} {\bf 30(3)} 283 (1998)





\end{thebibliography}
\end{document}